\documentstyle[12pt]{article}

\def\R{\bf R}
\def\C{\bf C}
\def\Z{\bf Z}
\def\P{\bf P}
\begin{document}
\begin{center}
{\Large {\bf Positively Curved Complete Noncompact K\" ahler Manifolds}} 
\vskip 10mm Bing--Long Chen , Xi--Ping Zhu

\smallskip 
Department of Mathematics, Zhongshan University,\\Guangzhou 510275, P. R.
China, and\\ \smallskip 
The Institute of Mathematical Sciences, The Chinese University of Hong Kong,
Hong Kong\vskip 18mm {\large {\bf Abstract}}
\end{center}

In this paper we give a partial affirmative answer to a conjecture of Greene-Wu and Yau. We prove that a complete noncompact K\"ahler surface with positive and bounded sectional curvature and with finite analytic Chern number $c_{1}(M)^{2}$ is biholomorphic to ${\C}^2$. \vskip 20mm \baselineskip=20pt

\section*{\S1. Introduction}

\setcounter{section}{1} \setcounter{equation}{0}

\qquad The celebrated theorem of Cheeger--Gromoll--Meyer \cite{CG}, \cite{GM}
states that a complete noncompact Riemannian manifold with positive
sectional curvature is diffeomorphic to the Euclidean space. It is well-known that there is
a vast variety of biholomorphically distinct complex structures on {\R}$%
^{2n} $ for $n>1$ ( see \cite{BSW}, \cite{F} ). To understand the
relationship between the Riemannian structure and the complex structure on
complete noncompact manifolds, we restrict attention to K\"ahler manifolds
which has the effect of insuring a closer relationship between these two
structures. In \cite{GW1}, Greene and Wu proved that a complete noncompact
K\"ahler manifold with positive sectional curvature is Stein. This fact thus
motivated the following conjecture:\vskip 3mm{\bf \underline{Conjecture}} (
Greene--Wu \cite{GW2} and Yau \cite{Y1} )\quad A complete noncompact K\"ahler manifold of
positive sectional curvature is biholomorphic to a complex Euclidean space. 
\vskip 3mm In \cite{Mo1}, among other things, Mok gave the first partial
affirmative answer to the conjecture for complex two--dimensional manifolds
with maximal volume growth.\vskip 3mm{\bf \underline{Theorem}} ( Mok \cite
{Mo1} )\quad Suppose $M$ is a complete noncompact K\"ahler manifold of
complex dimension $n=2$. Suppose also $M$ has positive sectional curvature
and satisfies\vskip 0.5mm 
$$
\begin{array}{ll}
\bigbreak (i)\qquad & 0<scalar\ \ curvature\leq 
\displaystyle \frac C{d^2(x_0,x)}\ , \\ \bigbreak (ii)\qquad & 
Volume(B(x_0,r))\geq cr^{2n}\ ,\qquad 0\leq r<+\infty \ , 
\end{array}
$$
where $B(x_0,r)$ and $d(x_0,x)$ denote respectively geodesic balls and
geodesic distances, c, C are some positive constants. Then $M$ is
biholomorphic to {\C}$^2$.\vskip 3mm Denote by $Ric$ the Ricci curvature
form of $M$. As remarked in Mok \cite{Mo3}, the conditions $(i)$ and $(ii)$
imply that the integral $\int_M{Ric^n}$, the analytic Chern number $%
c_{1}(M)^n$, is finite. 

In his paper \cite{T}, To gave a
generalization of the above result to nonmaximal volume growth manifolds. More precisely, it was proved in \cite{T} that if $M$ is a
complete noncompact K\"ahler manifold of positive sectional curvature with
complex dimension $n=2$ and suppose for some base point $x_0\in M$ that there
exist positive constants $C_1,\ C_2$ and $p$ such that\vskip 0.5mm 
$$
\begin{array}{ll}
\bigbreak (i)'\qquad & 0<scalar\ \ curvature\leq 
\displaystyle \frac{C_1}{d^p(x_0,x)}\ ,  \\ 
(ii)'\qquad & \displaystyle \int\limits_{B(x_0,r)}\frac 1{\left(
1+d(x_0,x)\right) ^{np}}dx\leq C_2\log (r+2)\ ,\qquad 0\leq r<+\infty \,
 \\ \bigbreak (iii) & 
\displaystyle \int_MRic^n<+\infty , 
\end{array}
$$
then $M$ is biholomorphic to {\C}$^2$. Moreover in case of $p\geq 2$, the result is
valid without assuming condition $(ii)'$. Also from Yau's theorem (cf.\cite{ScY}) that complete noncompact manifolds with nonnegative Ricci curvature have at least linear volume growth we see that the positive constant $p$ need to be assumed not less that $\frac 1 n$.\vskip 3mm

It is likely that the assumption $(iii)$ is automatically satiafied for complete K\"ahler surfaces with positive sectional curvature. The reason is that on a complete noncompact real four-manifold with positive sectional curvature we have the generalized Cohn-Vossen inequality
$$\int_M \Theta \leq \chi ({\R}^4)< +\infty$$
where $\Theta$ is the Gauss-Bonnet-Chern integrand. It is well-known that the integrand $\Theta$ is positive and it seems that the exterior product $Ric^2$ is more or less comparable with the Gauss-Bonnet-Chern integrand $\Theta$. Meanwhile in views of Demailly's holomorphic Morse inequality \cite{D1} and the $L^2$-Riemann-Roch inequality in Nadel-Tsuji\cite{NT} (see also Tian\cite{Ti}), the assumption $(iii)$ is a natural condition for a complete K\"ahler manifold to be a quasi-projective manifold. However the assumptions on the curvature decay and the volume growth are more problemtic since they demand the geometry of the K\"ahler manifold at infinity to be somewhat uniform. The main purpose of this paper is to show that the assumption $(iii)$ alone is sufficient to guarantee that the K\"ahler surface is biholomorphic to {\C}$^2$.
\vskip 3mm{\bf \underline{Main Theorem}}\quad Let $M$ be a complex $n$%
--dimensional complete noncompact K\"ahler manifold with positive and
bounded sectional curvature. Suppose that%
$$
\int_MRic^n<+\infty . 
$$
Then $M$ is biholomorphic to a quasi--projective variety, and in case of
dimenion $n=2$, $M$ is biholomorphic to {\C}$^2$.
\vskip 3mm We remark that
there is a more ambitious conjecture due to Yau \cite{Y1}, \cite{Y2}, i.e.,
the question is to demonstrate that every complete noncompact K\"ahler
manifold with positive holomorphic bisectional curvature is biholomorphic to
the complex Euclidean space. But for such a K\"ahler manifold with positive
holomorphic bisectional curvature, one doesn't even know whether the
manifold is simply connected. Moreover it is also unknown whether the
K\"ahler manifold is Stein, which is a conjecture of Siu \cite{Siu}. In the
companion paper \cite{CTZ}, the authors and S. H. Tang gave a partial
affirmative answer to the above Yau's conjecture. We proved that given a
complete noncompact complex two--dimensional K\"ahler manifold $M$ of
positive and bounded holomorphic bisectional curvature, suppose its geodesic
balls have maximal volume and its scalar curvature decays to zero at
infinity in the average sense, then $M$ is biholomorphic to {\C}$^2$.

The basic idea to approach these conjectures is to compactify the manifold $%
M $ as a quasi--projective variety. Siu, Yau \cite{SiY} and Mok \cite{Mo1}
initiated this program by first using the $L^2$--method of
Andreotti--Vessentini and H{\"ormander to establish a Siegel's theorem for a
field of meromorphic functions and then desingularizing the ``birational''
map associated to the Siegel theorem. The crux is how to choose a
subfield of meromorphic functions with suitable estimates for the
desingularization. In \cite{Mo1} Mok solved the Poincar\'e--Lelong equation
to show that there is plenty of holomorphic functions of polynomial growth
and that the subfield of meromorphic functions arising from quotients of
holomorphic functions of polynomial growth has the desirable estimates. In 
\cite{CTZ} the Ricci flow was used to understand the topology of the
manifold and to deduce that this field generated by holomorphic functions of
polynomial growth still has the desirable estimates. In \cite{Mo3} and \cite
{T}, the quotient fields were defined from holomorphic plurianticanonical
sections with polynomial growth or satisfying certain integrability
conditions associated to the assumptions on the curvature decay and the
volume growth respectively.}

Note that one can construct $L^2$ holomorphic
sections of the plurianticanonical system $\left\{ K^{-q},\ q>0\right\} $
without any assumption on the volume growth and the curvature decay.
We observe in the present paper that the space of $L^2$ sections forms a graded
algebra. Moreover by using some techniques and estimates developed from the Ricci flow, we are able to derive desirable B\'ezout estimates for the intersections of the zero
divisors of holomorphic $L^2$ sections. Based on these estimates we can
bound the Gauss--Bonnet integrals of curves which are the intersection of
zero divisors of $L^2$ holomorphic plurianticanonical sections. The main
theorem will then be proved in the following way: we first use the B\'ezout estimates to
obtain a Siegel type theorem for the field of meromorphic functions of $M$
arising from the quotients of the $L^2$ sections; we then get a
``meromorphic'' map of $M$ into a projective algebraic variety; we next use
the bounds on the Gauss--Bonnet integrals of the ``algebraic" curves to show
that the meromorphic map is almost surjective; finally we desingularize the
map into a biholomorphism from $M$ to a quasi--projective variety. 

      The composition of this paper is as follows. In
Section 2 we collect some basic results and prove the space of $L^2$
holomorphic plurianticanonical sections forms a {{\Z}$^{+}$}--graded
algebra. In Section 3 we derive the B\'ezout estimates and obtain a 
gradient estimate by using the Ricci flow, and then we bound the Gauss-Bonnet integrals of the curves cut by $L^2$ holomorphic plurianticanonical sections. Section 4 is devoted to the proof of the main theorem. 

The authors are grateful to Professor S. T. Yau for his interest and encouragement. The second author would like to thank Professor G. Tian for helpful discussion. This work was partially supported by The IMS of The Chinese University of Hong Kong and the Foundation for Outstanding Young Scholars of China.

\section*{\S2. $L^2$ holomorphic plurianticanonical sections}

\setcounter{section}{2} \setcounter{equation}{0}

\qquad First of all, let us recall the standard $L^2$--estimates of $
\overline{\partial }$ of Andreotti--Vesentini \cite{AV} and H{\"ormander 
\cite{Ho} in the case of Hermitian holomorphic line bundles over }K\"ahler
manifolds.\vskip 3mm{\bf \underline{Theorem 2.1}} ( Andreotti--Vesentini 
\cite{AV}, H{\"ormander \cite{Ho}} )

Let $(M,\omega )$ be a complete K\"ahler manifold equipped with a K\"ahler
form $\omega $. Let $L$ be a Hermitian holomorphic line bundle on $M$ and
denote by $C(L)$ the curvature form of $L$. Let $\varphi $ be a smooth
function and $c(x)$ be a positive continuous functions on $M$ such that%
$$
\sqrt{-1}\partial \overline{\partial }\varphi +C(L)+Ric\geq c(x)\omega \ . 
$$
Suppose $f$ is a $\overline{\partial }$--closed smooth $L$--valued (0,1)
form such that%
$$
\int_M\frac{\left\| f\right\| ^2}ce^{-\varphi }<+\infty \ . 
$$
Then the equation $\overline{\partial }u=f$ has a smooth solution with the
following estimate%
$$
\int_M\left\| u\right\| ^2e^{-\varphi }\leq \int_M\frac{\left\| f\right\| ^2}%
ce^{-\varphi }\ . 
$$
\hfill$\Box $ \vskip 3mm We will need the following sub--mean--value
inequality which can be found in \cite{LT}. Actually the sub--mean--value
inequality can be obtained directly by using the Green formula and the
standard estimate of the Green function on a Riemannian manifold with
nonnegative Ricci curvature.\vskip 3mm{\bf \underline{Lemma 2.2}}\quad Let $%
M $ be a complete Riemannian manifold with nonnegative Ricci curvature.
Suppose $f$ is a nonnegative smooth subharmonic function. Then there exists a
constant $C$, depending only on the dimension, such that for any $x_0\in M$,
we have%
$$
f(x_0)\leq C\cdot \frac
1{Vol\left( B\left( x_0,a\right) \right) }\int_{B\left( x_0,a\right) }f(x)\
, 
$$
for all $a>0.$\hfill$\Box $\vskip 3mm Now let $M$ be a complete K\"ahler
manifold of complex dimension $n$ with positive Ricci curvature. Denote by $%
K $ the canonical line bundle and $\Gamma ^2\left( M,K^{-q}\right) $ the
space of square--integrable holomorphic sections of the plurianticanonical
line bundle $K^{-q}$. Here $q$ is a positive integer. Fix a point $x\in M$,
suppose $\left\{ z_1,\cdots ,z_n\right\} $ is a holomorphic coordinate
system at $x$ with $z_1(x)=\cdots =z_n(x)=0$. Let $\left\{ s_0^{\prime
},\cdots ,s_n^{\prime }\right\} $ be a system of local holomorphic sections
of $K^{-q}$ at $x$ with $s_0^{\prime }(x)\neq 0,\ s_1^{\prime }(x)=\cdots
=s_n^{\prime }(x)=0$ and $d\left( {s_i^{\prime }}/{s_0^{\prime }}\right)
=dz_i,\ (\ i=1,\cdots ,n\ ),$ near $x$. Without loss of generality, we may
assume $\left\{ \sum\limits_{i=1}^n\left| z_i\right| ^2<1\right\} $ be a
holomorphic coordinate ball contained in $M$. Let $\rho $ be a smooth
cut--off function on $M$ satisfying $\rho \equiv 1$ on $B\left( 0,\frac
12\right) $ and $\rho \equiv 0$ outside $B(0,1)$. Since the Ricci curvature
is strictly positive everywhere, one can find a positive integer $q$ such
that%
$$
(q+1)Ric+(n+1)\sqrt{-1}\partial \overline{\partial }\left( \rho \log
\sum\limits_{i=1}^n\left| z_i\right| ^2\right) >0\ . 
$$
Choose $L=K^{-q}$, $\varphi =(n+1)\rho \log \sum\limits_{i=1}^n\left|
z_i\right| ^2$. By applying Theorem 2.1, we know that there exist global
holomorphic sections $u_0,u_1,\cdots,u_n$ of $K^{-q}$ such that%
$$
\overline{\partial }u_i=\overline{\partial }\left( \rho s_i^{\prime }\right) 
$$
and%
$$
\int_M\left\| u_i\right\| ^2e^{-(n+1)\rho \log \sum\limits_{i=1}^n\left|
z_i\right| ^2}<+\infty \ , 
$$
for $i=0,1,\cdots,n\ .$

Set%
$$
s_i=\rho s_i^{\prime }-u_i\ ,\qquad i=0,1,\cdots ,n\ . 
$$
Then each $s_i\in \Gamma ^2\left( M,K^{-q}\right) $ and $s_0(x)\neq 0,\
d(s_i/s_0)=dz_i\ (\ i=1,\cdots ,n\ )$ near $x$. This says that the
meromorphic functions $s_1/s_0,\cdots ,s_n/s_0$ give a local holomorphic
coordinate system at $x.$

Similarly for arbitrary two points $x$ and $y$ in $M$, one can choose $s_0,
s_1\in \Gamma ^2\left( M,K^{-q}\right) $ with $s_0(x)\neq 0,\ s_0(y)\neq 0,\
s_1(x)\neq 0,\ $and $s_1(y)=0$ so that the meromorphic function $f=s_1/s_0$
is holomorphic at $x$ and $y$ and satisfies $f(x)\neq 0$ and $f(y)=0.$

Hence we have showed the space $\bigcup\limits_{q>0}\Gamma ^2\left(
M,K^{-q}\right) $ gives local holomorphic coordinates and separates points
on $M$. Furthermore, we will show that $\bigcup\limits_{q>0}\Gamma ^2\left(
M,K^{-q}\right) $ forms a graded algebra. More precisely, we have\vskip 3mm%
{\bf \underline{Proposition 2.3}}\quad Suppose $M$ has positive and bounded
sectional curvature. Then the space $\bigcup\limits_{q>0}\Gamma ^2\left(
M,K^{-q}\right) $ forms a {{\Z}$^{+}$}--graded algebra. Moreover if $s\in
\bigcup\limits_{q>0}\Gamma ^2\left( M,K^{-q}\right) $, then $\left\|
s\right\| $ is bounded and $\left\| \nabla s\right\| $ is square--integrable.%
\vskip 3mm{\bf \underline{Proof.}}\quad From the standard Bochner--Kodaira
formula, 
\begin{equation}
\label{2.1}\bigtriangleup \left\| s\right\| ^2=\left\| \nabla s\right\|
^2-qR\left\| s\right\| ^2\ , 
\end{equation}
where $R$ is the scalar curvature of $M.$

Suppose that the sectional curvature is bounded from above by a positive
constant $K_0$. It is easy to see from (\ref{2.1})%
$$
\bigtriangleup \left\| s\right\| ^2\geq -qn^2K_0\left\| s\right\| ^2.
$$
Let $G=e^{\sqrt{qn^2K_0}%
\tau } \left\| s\right\|^2 $ be a
function defined on the product manifold $\widetilde{M}=M\times {\R}$,
equipped with the product metric. Let $\widetilde{\bigtriangleup }%
=\bigtriangleup +\frac{\partial ^2}{\partial \tau ^2}$ be the
Laplacian operator. It is clear that $\widetilde{M}$ has nonnegative Ricci
curvature and we have%
$$\widetilde{\bigtriangleup }G\geq 0. $$
Then by Lemma 2.2 and the standard volume comparison we have for any $x_0\in M,$%
\begin{equation}
\label{2.3}\left\| s\right\| ^2(x_0)\leq C\cdot \frac{e^{\sqrt{qn^2K_0}\cdot
a}-e^{-\sqrt{qn^2K_0}\cdot a}}{\sqrt{qn^2K_0}\cdot a}\cdot \frac 1{Vol\left(
B(x_0,a)\right) }\int_{B(x_0,a)}\left\| s\right\| ^2\ , 
\end{equation}
for all $a>0$. Since $M$ has positive sectional curvature, it is well--known
from Gromoll and Meyer \cite{GM} that the injectivity radius of $M$ at $x_0$
satisfies the following estimate%
$$
inj\left( M,x_0\right) \geq \frac \pi {\sqrt{K_0}}\ . 
$$
In particular we have the volume estimate%
$$
Vol\left( B\left( x_0,\frac \pi {2\sqrt{K_0}}\right) \right) \geq C\left(
\frac \pi {2\sqrt{K_0}}\right) ^{2n} 
$$
for some positive constant $C$ depending only on the dimension. So letting $%
a=\frac \pi {2\sqrt{K_0}}$ in (\ref{2.3}) we conclude%
$$
\sup \limits_{x_0\in M}\left\| s\right\| ^2(x_0)\leq C\cdot \frac{e^{\sqrt{
\frac{qn^2}4}\cdot \pi }-e^{-\sqrt{\frac{qn^2}4}\cdot \pi }}{\sqrt{\frac{qn^2%
}4}\cdot \pi }\cdot \left( \frac{2\sqrt{K_0}}\pi \right) ^{2n}\int_M\left\|
s\right\| ^2<+\infty \ . 
$$
This shows that $\left\| s\right\| ^2$ is bounded on $M$. Thus the space $%
\bigcup\limits_{q>0}\Gamma ^2\left( M,K^{-q}\right) $ forms a {{\Z}$^{+}$}%
--graded algebra under standard addition and multiplication over the complex
numbers {\C}.

We now consider the $L^2$ estimate for the gradient of a holomorphic section 
$s\in \Gamma ^2\left( M,K^{-q}\right) .$ Suppose $r$ is the distance
function from a fixed point $x_0$ on $M$. It follows from (\ref{2.1}) that
for any $a>0,$%
\begin{eqnarray}
\int\limits_{B(x_0,a)}\left\| \nabla s\right\| ^2\left( 1-\frac ra\right) ^2
& = & \int\limits_{B(x_0,a)}\bigtriangleup \left\| s\right\| ^2\left(
1-\frac ra\right) ^2+\int\limits_{B(x_0,a)}qR\left\| s\right\| ^2\left(
1-\frac ra\right) ^2\nonumber \\  
& \leq  & \frac 4a\int\limits_{B(x_0,a)}\left\| \nabla s\right\| \cdot
\left\| s\right\| \cdot \left( 1-\frac ra\right)
+qn^2K_0\int\limits_M\left\| s\right\| ^2\nonumber \\  
& \leq  & \frac 12\int\limits_{B(x_0,a)}\left\| \nabla s\right\| ^2\left(
1-\frac ra\right) ^2+\left( \frac 8{a^2}+qn^2K_0\right) \int\limits_M\left\|
s\right\| ^2_.\nonumber
\end{eqnarray}
Letting $a\rightarrow +\infty $, we get 
\begin{equation}
\label{2.4}\int\limits_M\left\| \nabla s\right\| ^2\leq
2qn^2K_0\int\limits_M\left\| s\right\| ^2<+\infty \ . 
\end{equation}

\hfill$\Box $

\section*{\S 3. B\'ezout estimates and gradient estimate}

\setcounter{section}{3} \setcounter{equation}{0} \qquad Let us first recall
a cut--off function constructed in the book of Schoen and Yau \cite{ScY} (
see Theorem 1.4.2 in \cite{ScY} ).\vskip 3mm{\bf \underline{Lemma 3.1}} (
Schoen--Yau \cite{ScY}, see also Shi{\ \cite{Sh}} )

Suppose $M$ is a complete Riemannian manifold with nonnegative and bounded
sectional curvature. Then for any $x_0\in M,\ a>0$, there exists a positive
function $\varphi _a(x)\in C^\infty (M)$ such that%
\begin{eqnarray}
(\alpha ) &  & \exp \left( -C\left( 1+\frac{d(x,x_0)}a\right) \right) \leq \varphi _a(x)\leq \exp \left( -\left( 1+\frac{d(x,x_0)}a\right) \right) \ ,\nonumber  \\
(\beta ) &  & \left| \nabla \varphi _a(x)\right| \leq \frac Ca\cdot\varphi _a(x)\ ,\nonumber  \\
(\gamma ) &  & \left| \nabla _i\nabla _j\varphi _a(x)\right| \leq \frac C{a^2}\cdot\varphi _a(x)\ ,\nonumber
\end{eqnarray}
for all $x\in M$, where $C$ is a positive constant depending only on the
dimension.

\hfill$\Box $\vskip 3mm We now state and prove the first B\'ezout estimate
for $\bigcup\limits_{q>0}\Gamma ^2\left( M,K^{-q}\right) .$\vskip 3mm{\bf 
\underline{Proposition 3.2}}\quad Let $(M,\omega )$ be assumed as in the
Main Theorem. Let $q_1,\cdots ,q_k$ be positive integers and let $s_i\in
\Gamma ^2\left( M,K^{-q_i}\right) ,\ i=1,2,\cdots ,k$. For a sequence $%
\varepsilon _i>0,\ i=1,\cdots ,k$, define%
$$
\zeta _{\varepsilon _i}(s_i)=\sqrt{-1}\partial \overline{\partial }\log \left(
\left\| s_i\right\| ^2+\varepsilon _i^2\right) +q_iRic\ ,\qquad i=1,\cdots ,k\
. 
$$
Then 
\begin{equation}
\label{3.1}\int_M\zeta _{\varepsilon _1}(s_1)\wedge \cdots \wedge \zeta
_{\varepsilon _k}(s_k)\wedge Ric^{n-k}\leq Cq_1\cdots q_k\int_MRic^n<+\infty \
, 
\end{equation}
where $C$ is a positive constant depending only on $n.$\vskip 3mm{\bf 
\underline{Proof.}}\quad We first consider the case $k=1$. For sufficiently
small $\delta >0$ and $a\geq 1$, we know $\left\{ \varphi _a>\delta \right\}
\subset M$ is a compact domain with boundary $\partial \left\{ \varphi
_a>\delta \right\} =\left\{ \varphi _a=\delta \right\} $. The Poincar\'e--Lelong equation gives, in the sense of (1,1) currents,
$$\sqrt{-1}\partial \overline{\partial }\log \left\| s\right\| ^2 = [s=0] - qRic$$
where, by abuse of notation, $[s=0]$ also denotes the (1,1) current defined by the divisor $[s=0]$ counting multiplicity. 
From the Poincar\'e--Lelong equation and the positivity of plurianticanonical line
bundles we see from an easy computation that each $\zeta _{\varepsilon _i}(s_i)$ is a closed,
nonnegative (1.1) form on $M$. By Lemma 3.1, we have
$$
\begin{array}{l}
\bigbreak\displaystyle\qquad \int_{\left\{ \varphi _a>\delta \right\}
}\left( \varphi _a-\delta \right) ^2\zeta _\varepsilon (s)\wedge Ric^{n-1} \\ 
\bigbreak\displaystyle=\int\limits_{\left\{ \varphi _a>\delta \right\}
}\left( \varphi _a-\delta \right) ^2\sqrt{-1}\partial \overline{\partial }%
\log \left( \left\| s\right\| ^2+\varepsilon ^2\right) \wedge
Ric^{n-1}+\int\limits_{\left\{ \varphi _a>\delta \right\} }q\left( \varphi
_a-\delta \right) ^2Ric^n \\ \bigbreak\displaystyle\leq \frac
C{a^2}\int_{\left\{ \varphi _a>\delta \right\} }\log \left( 1+\frac{\left\|
s\right\| ^2}{\varepsilon ^2}\right) \omega \wedge Ric^{n-1}+q\int_{\left\{
\varphi _a>\delta \right\} }Ric^n \\ \displaystyle\leq \frac{CK_0^{n-1}}{%
a^2\varepsilon ^2}\int_M\left\| s\right\| ^2+q\int_MRic^n 
\end{array}
$$
where $K_0$ is the upper bound of the sectional curvature of $M$. Here and
below we denote by $C$ various positive constants depending only on $n$.
Then letting $\delta \rightarrow 0$ and $a\rightarrow +\infty $ we deduce
that%
$$
\int_M\zeta _\varepsilon (s)\wedge Ric^{n-1}\leq Cq\int_MRic^n\ . 
$$

We next consider the case $k=2$. By integrating by parts, we get%
\begin{eqnarray}
\label{3.2}
&&\int_{\left\{ \varphi _a>\delta \right\} }\left( \varphi _a-\delta \right) ^3\zeta _{\varepsilon _1}(s_1)\wedge \zeta _{\varepsilon _2}(s_2)\wedge Ric^{n-2}\nonumber\\
&\leq&\int_{\left\{ \varphi _a>\delta \right\} }\left( \varphi _a-\delta \right) ^3\sqrt{-1}\partial \overline{\partial }\log \left( \left\| s_2\right\| ^2+\varepsilon _2^2\right) \wedge \zeta _{\varepsilon _1}(s_1)\wedge Ric^{n-2}\nonumber  \\
&&\qquad+\int_{\left\{ \varphi _a>\delta \right\} }q_2\zeta _{\varepsilon _1}(s_1)\wedge Ric^{n-1}\nonumber  \\
&=&\int_{\left\{ \varphi _a>\delta \right\} }\log \left( 1+\frac{\left\| s_2\right\| ^2}{\varepsilon _2^2}\right) \cdot \sqrt{-1}\partial \overline{\partial }\left( \left( \varphi _a-\delta \right) ^3\right) \wedge \zeta _{\varepsilon _1}(s_1)\wedge Ric^{n-2}\nonumber\\
&&\qquad+\int_{\left\{ \varphi _a>\delta \right\} }q_2\zeta _{\varepsilon _1}(s_1)\wedge Ric^{n-1}\ .
\end{eqnarray}
The second term of the RHS of (\ref{3.2}) is bounded from above by $%
Cq_1q_2\int_MRic^n$ by the previous estimate for the case $k=1$. To get the
bound on the first term of RHS of (\ref{3.2}), we now derive an integral
estimate for the gradient. Recall that for each $i,$%
$$
\zeta _{\varepsilon _i}(s_i)=\sqrt{-1}\partial \overline{\partial }\log \left(
\left\| s_i\right\| ^2+\varepsilon _i^2\right) +q_iRic\geq 0\ . 
$$
Thus we have%
$$
\int\limits_{\left\{ \varphi _a>\delta \right\} }\left( \varphi _a-\delta
\right) ^2\log \left( 1+\frac{\left\| s_i\right\| ^2}{\varepsilon _i^2}\right)
\left( \sqrt{-1}\partial \overline{\partial }\log \left( \left\| s_i\right\|
^2+\varepsilon _i^2\right) +q_iRic\right) \wedge \omega ^{n-1}\geq 0\ . 
$$
Integrating by parts and using Lemma 3.1, we get%
\begin{eqnarray}
&&\int_{\left\{ \varphi _a>\delta \right\} }\left( \varphi _a-\delta \right) ^2\sqrt{-1}\partial \log \left( 1+\frac{\left\| s_i\right\| ^2}
{\varepsilon _i^2}\right) \wedge \overline{\partial }\log \left( 1+\frac{\left\|s_i\right\|^2}{\varepsilon _i^2}\right)\wedge\omega^{n-1}
\nonumber\\
&\leq& \int_{\left\{ \varphi _a>\delta \right\} }\frac Ca\left( \varphi _a-\delta \right) \log \left( 1+\frac{\left\| s_i\right\| ^2}{\varepsilon _i^2}\right) \left\| \nabla \log \left( 1+\frac{\left\| s_i\right\| ^2}{\varepsilon _i^2}\right) \right\| \omega ^n\nonumber\\
&&\qquad+\int_{\left\{ \varphi _a>\delta \right\} }q_i\left( \varphi _a-\delta \right) ^2\log \left( 1+\frac{\left\| s_i\right\| ^2}{\varepsilon _i^2}\right) Ric\wedge \omega ^{n-1}\ .\nonumber
\end{eqnarray}
Recall from Proposition 2.3 that each $s_i$ is bounded on $M$. By using
Cauchy--Schwarz inequality, we get%
\begin{eqnarray}
&&\int_{\left\{ \varphi _a>\delta \right\} }\left( \varphi _a-\delta \right) ^2\left\| \nabla \log \left( 1+\frac{\left\| s_i\right\| ^2}{\varepsilon _i^2}\right) \right\| ^2\omega ^n\nonumber\\
&\leq&\frac C{a^2}\int_{\left\{ \varphi _a>\delta \right\} }\left( \log \left( 1+\frac{\left\| s_i\right\| ^2}{\varepsilon _i^2}\right) \right) ^2\omega ^n+Cq_iK_0\int_{\left\{ \varphi _a>\delta \right\} }\log \left( 1+\frac{\left\| s_i\right\| ^2}{\varepsilon _i^2}\right) \omega ^n\nonumber\\
&\leq&\frac C{a^2\varepsilon _i^4}\left( \sup \limits_M\left\| s_i\right\| ^2\right) \int_{\left\{ \varphi _a>\delta \right\} }\left\| 
s_i\right\| ^2\omega ^n+\frac{Cq_iK_0}{\varepsilon _i^2}\int_{\left\{ \varphi _a>\delta \right\} }\left\| s_i\right\| ^2\omega ^n\ .\nonumber
\end{eqnarray}
Thus letting $\delta \rightarrow 0$ and $a\rightarrow +\infty $ we deduce the integral estimate 
\begin{equation}
\label{3.3}\int_{M\ }\left\| \nabla \log \left( 1+\frac{\left\| s_i\right\|
^2}{\varepsilon _i^2}\right) \right\| ^2\omega ^n\leq \frac{Cq_iK_0}{\varepsilon
_i^2}\int_{M\ }\left\| s_i\right\| ^2\omega ^n\ . 
\end{equation}
Then the first term on the RHS of (\ref{3.2}) can be estimated as follows,%
\begin{eqnarray}\label{3.4}
&&\int_{\left\{ \varphi _a>\delta \right\} \ }\log \left( 1+\frac{\left\| s_2\right\| ^2}{\varepsilon _2^2}\right) \sqrt{-1}\partial
\overline{\partial }\left(\left( \varphi _a-\delta \right) ^3\right)\wedge \zeta _{\varepsilon _1}(s_1)\wedge Ric^{n-2}\nonumber\\
&\leq&\frac C{a^2}\int_{\left\{ \varphi _a>\delta \right\} \ }\left( \varphi _a-\delta \right)\log \left( 1+\frac{\left\| s_2\right\|^2}
{\varepsilon _2^2}\right) \cdot \sqrt{-1}\partial \overline{\partial }\log \left( \left\| s_1\right\| ^2+\varepsilon _1^2\right) \nonumber\\
&&\qquad\wedge \omega \wedge Ric^{n-2}+\frac{Cq_1}{a^2}\int_{\left\{ \varphi _a>\delta \right\} \ }\left( \varphi _a-\delta \right) \log \left( 1+\frac{\left\| s_2\right\| ^2}{\varepsilon _2^2}\right) \omega \wedge Ric^{n-1}\nonumber\\
&\leq&\frac{CK_0^{n-2}}{a^3}\int_{\left\{ \varphi _a>\delta \right\} \ }\log \left( 1+\frac{\left\| s_2\right\| ^2}{\varepsilon _2^2}
\right) \cdot \left\| \nabla \log \left( \left\| s_1\right\| ^2+\varepsilon _1^2\right) \right\| \omega ^n+\frac{CK_0^{n-2}}{a^2}\cdot
\nonumber\\
&&\qquad\int_{\left\{ \varphi _a>\delta \right\} \ }\left( \varphi _a-\delta \right) \left\| \nabla \log \left( \left\| s_1\right\| ^2+\varepsilon_1^2\right) \right\| \cdot \left\| \nabla \log \left( \left\| s_2\right\| ^2+\varepsilon _2^2\right) \right\| \omega ^n\nonumber\\
&&\qquad+\frac{Cq_1K_0^{n-1}}{a^2}\int_{\left\{ \varphi _a>\delta \right\} \ }\left( \varphi _a-\delta \right) \log \left( 1+\frac{\left\| s_2\right\| ^2}{\varepsilon _2^2}\right) \omega ^n\nonumber\\
&\leq&\frac{CK_0^{n-2}}{a^3}\int_{\left\{ \varphi _a>\delta \right\} \ }\left( \log \left( 1+\frac{\left\| s_2\right\| ^2}{\varepsilon _2^2}\right) \right) ^2\omega ^n+\left( \frac{CK_0^{n-2}}{a^3}+\frac{CK_0^{n-2}}{a^2}\right)\cdot\nonumber\\
&&\qquad\int_{\left\{ \varphi _a>\delta \right\} \ }\left( \left\| \nabla \log \left( \left\| s_1\right\| ^2+\varepsilon _1^2\right) \right\| ^2+\left\| \nabla \log \left( \left\| s_2\right\| ^2+\varepsilon _2^2\right) \right\| ^2\right) \omega ^n\nonumber\\
&&\qquad+\frac{Cq_1K_0^{n-1}}{a^2}\int_{\left\{ \varphi _a>\delta \right\} \ }\log \left( 1+\frac{\left\| s_2\right\| ^2}{\varepsilon _2^2}\right) \omega ^n\nonumber\\
&\leq&\frac{CK_0^{n-2}}{a^3\varepsilon _2^2}\sup \limits_M\left( \log \left( 1+\frac{\left\| s_2\right\| ^2}{\varepsilon _2^2}\right) \right) \int\limits_{\left\{ \varphi _a>\delta \right\} \ }\left\| s_2\right\| ^2\omega ^n+\frac{CK_0^{n-2}}{a^2}\left( 1+\frac 1a\right)\cdot\nonumber\\
&&\qquad\int_{\left\{ \varphi _a>\delta \right\} \ }\left(\left\| \nabla \log \left( 1+\frac{\left\| s_1\right\| ^2}{\varepsilon _1^2}\right) \right\|^2+\left\| \nabla \log \left( 1+\frac{\left\| s_2\right\| ^2}{\varepsilon _2^2}\right) \right\|^2\right)\omega ^n\nonumber\\
&&\qquad+\frac{Cq_1K_0^{n-1}}{a^2\varepsilon _2^2}\int_{\left\{ \varphi _a>\delta \right\} \ }\left\| s_2\right\| ^2\omega ^n\ .
\end{eqnarray}
Letting $\delta \rightarrow 0$ and $a\rightarrow +\infty $, we obtain from (%
\ref{3.2}), (\ref{3.3}), and (\ref{3.4}) that%
$$
\int_{M\ }\zeta _{\varepsilon _1}(s_1)\wedge \zeta _{\varepsilon _2}(s_2)\wedge
Ric^{n-2}\leq Cq_1q_2\int_MRic^n\ . 
$$

For the general case $k>2$, by inducting on $k$ and integrating by parts, we
have%
\begin{eqnarray}\label{3.5}
&&\int_{\left\{ \varphi _a>\delta \right\} \ }\left( \varphi _a-\delta \right) ^{k+1}\zeta _{\varepsilon _1}(s_1)\wedge \cdots \wedge \zeta _{\varepsilon _k}(s_k)\wedge Ric^{n-k}\nonumber\\
&=&\int_{\left\{ \varphi _a>\delta \right\} \ }\left( \varphi _a-\delta \right) ^{k+1}\sqrt{-1}\partial \overline{\partial }\log \left(1+
\frac{\left\|s_k\right\|^2}{\varepsilon _k^2}\right) \wedge\zeta _{\varepsilon _1}(s_1)\wedge \cdots \wedge \zeta _{\varepsilon_{k-1}}
(s_{k-1})\nonumber\\
&&\wedge Ric^{n-k}+q_k\int_{\left\{ \varphi _a>\delta \right\} \ }\left( \varphi _a-\delta \right) ^{k+1}\zeta _{\varepsilon _1}(s_1)\wedge \cdots \wedge \zeta _{\varepsilon _{k-1}}(s_{k-1})\wedge Ric^{n-k+1}\nonumber\\
&\leq&\frac C{a^2}\int_{\left\{ \varphi _a>\delta \right\} \ }\left( \varphi _a-\delta \right) ^{k-1}\log \left( 1+\frac{\left\| s_k\right\| ^2}{\varepsilon _k^2}\right) \zeta _{\varepsilon _1}(s_1)\wedge \cdots \wedge \zeta _{\varepsilon _{k-1}}(s_{k-1})\wedge\nonumber\\
&&Ric^{n-k}+q_k\int_{\left\{ \varphi _a>\delta \right\} \ }\left( \varphi _a-\delta \right) ^{k+1}\zeta _{\varepsilon _1}(s_1)\wedge \cdots \wedge \zeta _{\varepsilon _{k-1}}(s_{k-1})\wedge Ric^{n-k+1}\nonumber\\
&\leq&\frac C{a^2}\int_{\left\{ \varphi _a>\delta \right\} \ }\left( \varphi _a-\delta \right) ^{k-1}\log \left( 1+\frac{\left\| s_k\right\| ^2}{\varepsilon _k^2}\right) \sqrt{-1}\partial \overline{\partial }\log \left( 1+\frac{\left\| s_{k-1}\right\| ^2}{\varepsilon _{k-1}^2}\right) \wedge\nonumber\\
&&\qquad\qquad\zeta _{\varepsilon _1}(s_1)\wedge \cdots \wedge \zeta _{\varepsilon _{k-2}}(s_{k-2})\wedge\omega\wedge Ric^{n-k}\nonumber\\
&&\qquad+\frac{Cq_{k-1}}{a^2}\int_{\left\{ \varphi _a>\delta \right\} \ }\left( \varphi _a-\delta \right) ^{k-1}\log \left( 1+\frac{\left\| s_k\right\| ^2}{\varepsilon _k^2}\right) \zeta _{\varepsilon _1}(s_1)\wedge\cdots \wedge\nonumber\\
&&\qquad\qquad\zeta _{\varepsilon _{k-2}}(s_{k-2})\wedge \omega \wedge Ric^{n-k+1}+Cq_1\cdots q_k\int_MRic^n\nonumber\\
&\leq&\frac{\widetilde{C}}{a^2}\int_{\left\{ \varphi _a>\delta \right\} \ }\left( \varphi _a-\delta \right) ^{k-2}\cdot\nonumber\\
&&\left( \left\| \nabla \log \left( 1+\frac{\left\| s_{k-1}\right\| ^2}{\varepsilon _{k-1}^2}\right) \right\| ^2+\left\| \nabla \log \left( 1+\frac{\left\| s_k\right\| ^2}{\varepsilon _k^2}\right) \right\| ^2+\log \left( 1+\frac{\left\| s_k\right\| ^2}{\varepsilon _k^2}\right) \right)\nonumber\\
&&\qquad\wedge \zeta _{\varepsilon _1}(s_1)\wedge \cdots \wedge \zeta _{\varepsilon _{k-2}}(s_{k-2})\wedge\omega^{n-k+2}+Cq_1\cdots q_k\int_MRic^n\ .
\end{eqnarray}
Here and below we denote by $\widetilde{C}$ various positive constants
depending only on $n,\varepsilon _1,\cdots ,\varepsilon _k,q_1,\cdots ,q_k,\sup
\limits_M\left\| s_1\right\| ^2,\cdots ,\sup \limits_M\left\| s_k\right\| ^2$
and $K_0.$

We now show that the following two estimates 
\begin{equation}
\label{3.6}\int_{M\ }\log \left( 1+\frac{\left\| s_k\right\| ^2}{\varepsilon
_k^2}\right) \zeta _{\varepsilon _1}(s_1)\wedge \cdots \wedge \zeta _{\varepsilon
_{k-2}}(s_{k-2})\wedge \omega ^{n-k+2}<+\infty 
\end{equation}
and 
\begin{equation}
\label{3.7}\int_{M\ }\left\| \nabla \log \left( 1+\frac{\left\|
s_{k-1}\right\| ^2}{\varepsilon _{k-1}^2}\right) \right\| ^2\zeta _{\varepsilon
_1}(s_1)\wedge \cdots \wedge \zeta _{\varepsilon _{k-2}}(s_{k-2})\wedge \omega
^{n-k+2}<+\infty 
\end{equation}
hold by induction. The estimate (\ref{3.3}) tells us that these integration
are finite for the case $k=2$. Since%
\begin{eqnarray}
&&\int_{\left\{ \varphi _a>\delta \right\} \ }\left( \varphi _a-\delta \right) ^2\log \left( 1+\frac{\left\| s_k\right\| ^2}{\varepsilon _k^2}\right) \zeta _{\varepsilon _1}(s_1)\wedge \cdots \wedge \zeta _{\varepsilon _{k-2}}(s_{k-2})\wedge \omega ^{n-k+2}\nonumber\\
&\leq&\frac Ca\int_{\left\{ \varphi _a>\delta \right\} \ }\left( \varphi _a-\delta \right)\log \left( 1+\frac{\left\| s_k\right\| ^2}{\varepsilon _k^2}\right) \left\| \nabla \log \left( 1+\frac{\left\| s_{k-2}\right\| ^2}{\varepsilon _{k-2}^2}\right) \right\|\nonumber\\
&&\qquad\qquad\zeta _{\varepsilon _1}(s_1)\wedge \cdots \wedge \zeta _{\varepsilon _{k-3}}(s_{k-3})\wedge\omega^{n-k+3}\nonumber\\
&&+\int_{\left\{ \varphi _a>\delta \right\} }\left( \varphi _a-\delta \right) ^2\left\| \nabla \log \left( 1+\frac{\left\| s_k\right\| ^2}{\varepsilon _k^2}\right) \right\| \cdot \left\| \nabla \log \left( 1+\frac{\left\| s_{k-2}\right\| ^2}{\varepsilon _{k-2}^2}\right) \right\|\nonumber\\
&&\qquad\qquad\zeta _{\varepsilon _1}(s_1)\wedge \cdots \wedge \zeta _{\varepsilon _{k-3}}(s_{k-3})\wedge \omega ^{n-k+3}\nonumber\\
&&+q_{k-2}K_0\int_{\left\{ \varphi _a>\delta \right\} \ }\left( \varphi _a-\delta \right) ^2\log \left( 1+\frac{\left\| s_k\right\| ^2}{\varepsilon _k^2}\right) \zeta _{\varepsilon _1}(s_1)\wedge \cdots \wedge\nonumber\\
&&\qquad\qquad\zeta _{\varepsilon _{k-3}}(s_{k-3})\wedge\omega^{n-k+3}\ ,\nonumber
\end{eqnarray}
it follows (\ref{3.6}) directly by induction on the both estimates. To get (%
\ref{3.7}), we recall that%
$$
\zeta _{\varepsilon _{k-1}}(s_{k-1})=\sqrt{-1}\partial \overline{\partial }\log
\left( \left\| s_{k-1}\right\| ^2+\varepsilon _{k-1}^2\right) +q_{k-1}Ric\geq
0\ . 
$$
Thus%
\begin{eqnarray}
\int_{\left\{ \varphi _a>\delta \right\} \ }\left( \varphi _a-\delta \right)
^2\log \left( 1+\frac{\left\| s_{k-1}\right\| ^2}{\varepsilon _{k-1}^2}%
\right) ( \sqrt{-1}\partial \overline{\partial }\log \left( \left\|
s_{k-1}\right\| ^2+\varepsilon _{k-1}^2\right)&&\nonumber\\
+q_{k-1}Ric) \wedge\zeta _{\varepsilon _1}(s_1)\wedge \cdots \wedge \zeta _{\varepsilon
_{k-2}}(s_{k-2})\wedge \omega ^{n-k+1}&\geq&0\ .\nonumber
\end{eqnarray}
Integrating by parts we have%
\begin{eqnarray}
&&\int_{\left\{ \varphi _a>\delta \right\} \ }\left( \varphi _a-\delta \right) ^2\sqrt{-1}\partial \log \left( 1+\frac{\left\| s_{k-1}\right\| ^2}{\varepsilon _{k-1}^2}\right) \wedge \overline{\partial }\log \left( 1+\frac{\left\| s_{k-1}\right\| ^2}{\varepsilon _{k-1}^2}\right) \wedge\nonumber\\
&&\qquad\qquad\zeta _{\varepsilon _1}(s_1)\wedge \cdots \wedge \zeta _{\varepsilon _{k-2}}(s_{k-2})\wedge \omega ^{n-k+1}\nonumber\\
&\leq&\frac Ca\int_{\left\{ \varphi _a>\delta \right\} \ }\left( \varphi _a-\delta \right) \log \left( 1+\frac{\left\| s_{k-1}\right\| ^2}{\varepsilon _{k-1}^2}\right) \left\| \nabla \log \left( 1+\frac{\left\| s_{k-1}\right\| ^2}{\varepsilon _{k-1}^2}\right) \right\| \zeta _{\varepsilon _1}(s_1)\wedge\nonumber\\
&&\qquad\cdots \wedge\zeta _{\varepsilon _{k-2}}(s_{k-2})\wedge \omega ^{n-k+2}+q_{k-1}K_0\cdot\nonumber\\
&&\int\limits_{\left\{ \varphi _a>\delta \right\} \ }\left( \varphi _a-\delta \right) ^2\log \left( 1+\frac{\left\| s_{k-1}\right\| ^2}{\varepsilon _{k-1}^2}\right) \zeta _{\varepsilon _1}(s_1)\wedge \cdots \wedge\zeta _{\varepsilon _{k-2}}(s_{k-2})\wedge \omega ^{n-k+2}\ .\nonumber
\end{eqnarray}
Applying Cauchy--Schwarz inequality, we get 
$$
\begin{array}{l}
\bigbreak\displaystyle\int\limits_{\left\{ \varphi _a>\delta \right\} \
}\left( \varphi _a-\delta \right) ^2\left\| \nabla \log \left( 1+\frac{%
\left\| s_{k-1}\right\| ^2}{\varepsilon _{k-1}^2}\right) \right\| ^2\zeta
_{\varepsilon _1}(s_1)\wedge \cdots \wedge \zeta _{\varepsilon
_{k-2}}(s_{k-2})\wedge \omega ^{n-k+2} \\ \bigbreak\displaystyle\leq \frac
C{a^2}\int_{\left\{ \varphi _a>\delta \right\} \ }\left( \log \left( 1+\frac{%
\left\| s_{k-1}\right\| ^2}{\varepsilon _{k-1}^2}\right) \right) ^2\zeta
_{\varepsilon _1}(s_1)\wedge \cdots \wedge \zeta _{\varepsilon
_{k-2}}(s_{k-2})\wedge \omega ^{n-k+2} \\ \bigbreak\displaystyle\quad
+q_{k-1}K_0\int_{\left\{ \varphi _a>\delta \right\} \ }\left( \varphi
_a-\delta \right) ^2\log \left( 1+\frac{\left\| s_{k-1}\right\| ^2}{\varepsilon
_{k-1}^2}\right) \zeta _{\varepsilon _1}(s_1)\wedge \cdots \wedge \zeta
_{\varepsilon _{k-2}}(s_{k-2})\wedge \\ \bigbreak\displaystyle\quad \omega
^{n-k+2} \\ \displaystyle\leq \left( \frac{\widetilde{C}}{a^2}+\widetilde{C}%
\right) \int_{M\ }\log \left( 1+\frac{\left\| s_{k-1}\right\| ^2}{\varepsilon
_{k-1}^2}\right) \zeta _{\varepsilon _1}(s_1)\wedge \cdots \wedge \zeta
_{\varepsilon _{k-2}}(s_{k-2})\wedge \omega ^{n-k+2}\ . 
\end{array}
$$
Hence we obtain (\ref{3.7}) by using (\ref{3.6}).

We thus let $\delta \rightarrow 0$ and $a\rightarrow +\infty $ in (\ref{3.5}%
) and use (\ref{3.6}) and (\ref{3.7}) to obtain the desired estimate (\ref
{3.1}).\hfill$\Box $\vskip 3mm We will need to bound the Gauss--Bonnet
integrals of curves which are the intersection of zero divisors of sections
in $\bigcup\limits_{q>0}\Gamma ^2\left( M,K^{-q}\right) $. To this end we
want to establish the associated B\'ezout estimate over projectivized
tangent bundle.

We denote by $\pi :{\P}TM\rightarrow M$ the projectivied tangent bundle and $%
L$ the tautological line bundle. Fix a point $x_0\in M$, let $\left\{
z_1,\cdots,z_n\right\} $ be a holomorphic coordinate system at $x_0$ with $%
z_1(x_0)=\cdots =z_n(x_0)=0$. For any tangent vector $v=v_1\frac \partial
{\partial z_1}+\cdots +v_n\frac \partial {\partial z_n}$ with $v_n\neq 0$,
let $u_1=\frac{v_1}{v_n},\cdots ,u_{n-1}=\frac{v_{n-1}}{v_n}$. Then $\left\{
z_1,\cdots,z_n,u_1,\cdots ,u_{n-1}\right\} $ forms a local holomorphic
coordinate system on ${\P}TM$ at $\left( x_0,\left[ \frac \partial {\partial
z_1},\cdots ,\frac \partial {\partial z_n}\right] \right) $. Equip ${\P}TM$
and $L$ with the induced metric from $M$. Then a direct computation ( see 
\cite{Mo2} ) gives 
\begin{equation}
\label{3.8}c_1(L)=-\sqrt{-1}\sum\limits_{1\leq \gamma \leq n-1}du^\gamma
\wedge d\overline{u}^\gamma +\sqrt{-1}\sum\limits_{1\leq i,j\leq n}R_{n 
\overline{n}i\overline{j}}dz^i\wedge d\overline{z}^j 
\end{equation}
where $c_1(L)$ is the first Chern form of the tautological line bundle. Since
the sectional curvature of $M$ is positive, we thus have 
\begin{equation}
\label{3.9}-c_1(L)+2\cdot \left( \pi ^{*}Ric\right) \geq \sum\limits_{1\leq
\gamma \leq n-1}du^\gamma \wedge d\overline{u}^\gamma +\sum\limits_{1\leq
i,j\leq n}R_{i\overline{j}}dz^i\wedge d\overline{z}^j>0\ . 
\end{equation}
Hence $L^{*}\otimes \pi ^{*}K^{-q}$ is a positive line bundle over ${\P}TM $
with bounded curvature for $p\geq 2.$

Set%
$$
\nu =-c_1(L)+2\cdot \left( \pi ^{*}Ric\right) \ . 
$$
Now $\nu $ is a positive, closed (1.1) form on ${\P }TM$. Thus, $\nu $ can
be regarded as a ( possible incomplete ) K\"ahler metric on ${\P }TM.$

Recall that $\varphi _a$ is the cut--off function in Lemma 3.1. Now we
choose $\pi ^{*}\varphi _a$ as a cut--off function on ${\P }TM$. We remark that in the proof of Lemma 3.2 the main ingredients are the boundedness of the $L^2$ norm and $L^\infty$ norm of the section $s$. So one can
proceed exactly as in the proof of Proposition 3.2 to obtain the following
B\'ezout estimate.\vskip 3mm{\bf \underline{Proposition 3.3}}\quad Let ($%
M,\omega $) be assumed as in the Main Theorem. Let $s_i\in \Gamma ^2\left( {%
\P }TM,\pi ^{*}K^{-p_i}\right) ,\ i=1,\cdots ,k,\ p_1,\cdots ,p_k\geq 2$ and 
$t_j\in \Gamma ^2({\P }TM,L^{*}\otimes \\\pi ^{*}K^{-q_j}),\ j=1,\cdots ,l,\
q_1,\cdots ,q_l\geq 2$. Suppose also that all $\left\| s_i\right\| $ and $%
\left\| t_j\right\| $ are bounded. For any sequence $\varepsilon _i>0,\
i=1,\cdots ,k+l$, define%
$$
\zeta _{\varepsilon _i}(s_i)=\sqrt{-1}\partial \overline{\partial }\log \left(
\left\| s_i\right\| ^2+\varepsilon _i^2\right) +c_1\left( \pi
^{*}K^{-p_i}\right) \ ,\qquad i=1,\cdots ,k\ , 
$$
and%
$$
\eta _{\varepsilon _{j+k}}(t_j)=\sqrt{-1}\partial \overline{\partial }\log
\left( \left\| t_j\right\| ^2+\varepsilon _{j+k}^2\right) +c_1\left(
L^{*}\otimes \pi ^{*}K^{-q_j}\right) \ ,\quad j=1,\cdots ,l\ , 
$$
where $c_1\left( \pi ^{*}K^{-p_i}\right) $ and $c_1\left( L^{*}\otimes \pi
^{*}K^{-q_i}\right) $ are the first Chern forms of the corresponding line
bundles. Then%
\begin{eqnarray}
\label{3.10}&&\int_{{\P}TM}\zeta _{\varepsilon _1}(s_1)\wedge \cdots \wedge
\zeta _{\varepsilon _k}(s_k)\wedge \eta _{\varepsilon _{k+1}}(t_1)\wedge \cdots
\wedge \eta _{\varepsilon _{k+l}}(t_l)\wedge v^{2n-1-k-l}\nonumber\\
&\leq&C\ p_1\cdots
p_kq_1\cdots q_l\int_{{\P}TM}\nu^{2n-1}\ ,
\end{eqnarray}
where $C$ is a positive constant depending only on $n.$\hfill$\Box $

We remark that ${\P }TM$ has finite volume with respect to the K\"ahler
metric $\nu ,\ i.e.,$%
\begin{eqnarray}
\label{3.11}
\int_{{\P}TM}\nu^{2n-1}&\leq&Const\cdot \int_M\left( \int_{{\P}T_xM}\omega _{FS}\right) Ric^n\nonumber  \\
&\leq&Const\cdot Vol\left( {\P}^{n-1}\right) \cdot \int_MRic^n<+\infty \ ,
\end{eqnarray}
where we denote by $\omega _{FS}$ the Fubini--Study metric on ${\P }^{n-1}.$

We are ready to derive an estimate for the Gauss--Bonnet integrals of curves
which are the intersection of zero divisors of section in $%
\bigcup\limits_{q>0}\Gamma ^2\left( M,K^{-q}\right) $.

For any smooth holomorphic curve $S$ on $M$, let $\widehat{\theta }:S\rightarrow {%
\P }TM$ denote the lifting of $S$ to ${\P }TM$ defined by $\widehat{\theta }%
(x)=\left[ T_xS\right] $ for $x\in S$, where $\left[ T_xS\right] $ denotes
the element in ${\P }TM$ defined by $T_xS$. Let $\widehat{S}$ denote the
image $\widehat{\theta }(S)$. Suppose $S$ lies in the intersection of zero
divisors of $n-1$ holomorphic sections $t_1,\cdots ,t_{n-1}$, where $t_i\in
\Gamma ^2\left( M,K^{-q}\right) $ ( for some positive integer $q\geq 2$ ). Without loss of generality, we may assume that $q$ is sufficiently large. By $L^2$-estimate of $\overline{\partial}$-operator as in Section 2, we can choose an $L^2$ holomorphic plurianticanonical section $t_{0}$ of the same degree $q$ such that $t_{0}$ does not vanish identically on each zero divisor of holomorphic section $t_{i}, 1 \leq i \leq n-1$.
As shown in Section 4 of \cite{Mo2}, all $d\left( t_i/t_0\right) \cdot
t_0^2=t_0\nabla t_i-t_i\nabla t_0$ are holomorphic sections of $L^{*}\otimes
\pi ^{*}K^{-2q}$ over ${\P }TM$ for $i=1,\cdots ,n-1$, and $\widehat{S}$
lies in the intersection of $2n-2$ zero divisors of $\pi ^{*}t_1,\cdots ,\pi
^{*}t_{n-1},\ t_0\nabla t_1-t_1\nabla t_0,\cdots ,t_0\nabla
t_{n-1}-t_{n-1}\nabla t_0$. In order to apply Proposition 3.3, we will
show that all $\left\| \nabla t_i\right\| ,\ i=0,\cdots ,n-1$, are bounded
on $M$. Now we use the Ricci flow to get the following gradient estimate for 
$L^2$ holomorphic plurianticanonical sections, which has its own interest in analysis.
\vskip 3mm{\bf \underline{%
Proposition 3.4}}\quad Suppose $M$ has positive and bounded sectional
curvature. Let $s$ be a holomorphic section belonging $\Gamma ^2\left(
M,K^{-q}\right) $. Then $\left\| \nabla s\right\| $ is bounded on $M.$\vskip %
3mm{\bf \underline{Proof.}}\quad Let $\left( z_1,\cdots ,z_n\right) $ be a
local holomorphic coordinate system and write%
$$
s=f\left( \frac \partial {\partial z_1}\wedge \cdots \wedge \frac \partial
{\partial z_n}\right) ^q 
$$
locally for some holomorphic function $f$. Then $h=\det \left( g_{\alpha 
\overline{\beta }}\right) ^q$ is the Hermitian metric on $K^{-q}$ and $%
C(L)=- \sqrt{-1}\partial \overline{\partial }\log h=qRic$ is the curvature
form of $L=K^{-q}$ with respect to the coordinates $z_1,\cdots ,z_n$. The
covariant derivative of $s$ is given by 
\begin{equation}
\label{3.12}\nabla s=\left( \frac{\partial f}{\partial z_i}+qf\frac \partial
{\partial z_i}\log \det \left( g_{\alpha \overline{\beta }}\right) \right)
\otimes dz_i\otimes \left( \frac \partial {\partial z_i}\wedge \cdots \wedge
\frac \partial {\partial z_n}\right) ^q\ . 
\end{equation}
We use the Bochner trick to bound $\nabla s.$ Let $\left( z_1,\cdots
,z_n\right) $ be a normal holomorphic coordinate at a fixed point. We compute%
$$
\Vert \nabla s\Vert ^2=g^{i\overline{j}}\left( \frac{\partial f}{\partial z_i%
}+qf\frac \partial {\partial z_i}\log \det (g_{\alpha \overline{\beta }%
})\right) \overline{\left( \frac{\partial f}{\partial z_j}+qf\frac \partial
{\partial z_j}\log \det (g_{\alpha \overline{\beta }})\right) }\det
(g_{\alpha \overline{\beta }})_{,}^q 
$$
$$
\begin{array}{l}
\bigbreak\displaystyle\quad \frac \partial {\partial \overline{z}_l}\left\|
\nabla s\right\| ^2 \\ \bigbreak\displaystyle=\frac{\partial g^{i\overline{j}%
}}{\partial \overline{z}_l}\left( \frac{\partial f}{\partial z_i}+qf\frac
\partial {\partial z_i}\log \det \left( g_{\alpha \overline{\beta }}\right)
\right) \overline{\left( \frac{\partial f}{\partial z_j}+qf\frac \partial
{\partial z_j}\log \det \left( g_{\alpha \overline{\beta }}\right) \right) }%
\det \left( g_{\alpha \overline{\beta }}\right) ^q \\ \bigbreak\displaystyle%
+ \frac{\partial \det \left( g_{\alpha \overline{\beta }}\right) ^q}{%
\partial \overline{z}_l}g^{i\overline{j}}\left( \frac{\partial f}{\partial
z_i}+qf\frac \partial {\partial z_i}\log \det \left( g_{\alpha \overline{%
\beta }}\right) \right) \overline{\left( \frac{\partial f}{\partial z_j}%
+qf\frac \partial {\partial z_j}\log \det \left( g_{\alpha \overline{\beta }%
}\right) \right) } \\ \bigbreak\displaystyle+g^{i\overline{j}}\left( qf\frac{%
\partial ^2}{\partial z_i\partial \overline{z}_l}\log \det \left( g_{\alpha 
\overline{\beta }}\right) \right) \overline{\left( \frac{\partial f}{%
\partial z_j}+qf\frac \partial {\partial z_j}\log \det \left( g_{\alpha 
\overline{\beta }}\right) \right) }\det \left( g_{\alpha \overline{\beta }%
}\right) ^q \\ \bigbreak\displaystyle+g^{i\overline{j}}\left( \frac{\partial
f}{\partial z_i}+qf\frac \partial {\partial z_i}\log \det \left( g_{\alpha 
\overline{\beta }}\right) \right) \cdot \\ \bigbreak\displaystyle\quad 
\overline{\left( \frac{\partial ^2f}{\partial z_j\partial z_l}+qf\frac{%
\partial ^2}{\partial z_j\partial z_l}\log \det \left( g_{\alpha \overline{%
\beta }}\right) +q\frac{\partial f}{\partial z_l}\frac \partial {\partial
z_j}\log \det \left( g_{\alpha \overline{\beta }}\right) \right) }\det
\left( g_{\alpha \overline{\beta }}\right) ^q\ , 
\end{array}
$$
and 
$$
\begin{array}{l}
\bigbreak\displaystyle\quad \frac{\partial ^2}{\partial z_k\partial 
\overline{z}_l}\left\| \nabla s\right\| ^2 \\ \displaystyle=R_{k\overline{l}%
i \overline{j}}\left( \frac{\partial f}{\partial z_i}+qf\frac \partial
{\partial z_i}\log \det \left( g_{\alpha \overline{\beta }}\right) \right) 
\overline{\left( \frac{\partial f}{\partial z_j}+qf\frac \partial {\partial
z_j}\log \det \left( g_{\alpha \overline{\beta }}\right) \right) }\det
\left( g_{\alpha \overline{\beta }}\right) ^q 
\end{array}
$$
$$
+\frac{\partial ^2\det \left( g_{\alpha \overline{\beta }}\right) ^q}{%
\partial z_k\partial \overline{z}_l}\cdot g^{i\overline{j}}\left( \frac{%
\partial f}{\partial z_i}+qf\frac \partial {\partial z_i}\log \det \left(
g_{\alpha \overline{\beta }}\right) \right) \overline{\left( \frac{\partial
f }{\partial z_j}+qf\frac \partial {\partial z_j}\log \det \left( g_{\alpha 
\overline{\beta }}\right) \right) } 
$$
\begin{eqnarray}
&&+g^{i\overline{j}}\left( q\frac{\partial f}{\partial z_k}\frac{\partial ^2}{\partial z_i\partial \overline{z}_l}\log \det \left( g_{\alpha \overline{\beta }}\right) +qf\frac{\partial ^3}{\partial z_i\partial \overline{z}_l\partial z_k}\log \det \left( g_{\alpha \overline{\beta }}\right) \right)
\times\nonumber\\
&&\quad\overline{\left( \frac{\partial f}{\partial z_j}+qf\frac \partial {\partial
z_j}\log \det \left( g_{\alpha \overline{\beta }}\right) \right) }\det
\left( g_{\alpha \overline{\beta }}\right) ^q\nonumber\\ 
&&+g^{i\overline{j}}\left( qf\frac{\partial ^2}{\partial z_i\partial \overline{z}_l}\log \det \left( g_{\alpha \overline{\beta }}\right) \right) \overline{\left( qf\frac {\partial ^2}{\partial \overline{z}_k\partial z_j}\log \det
\left( g_{\alpha \overline{\beta }}\right) \right) }\det \left( g_{\alpha 
\overline{\beta }}\right) ^q\nonumber\\ 
&&+g^{i\overline{j}}\left( \frac{\partial ^2f}{\partial z_i\partial z_k}+q
\frac{\partial f}{\partial z_k}\frac \partial {\partial z_i}\log \det \left(
g_{\alpha \overline{\beta }}\right) +qf\frac{\partial ^2}{\partial
z_i\partial z_k}\log \det \left( g_{\alpha \overline{\beta }}\right) \right)
\times\nonumber\\ 
&&\quad\overline{\left( \frac{\partial ^2f}{\partial z_j\partial z_l}+qf\frac{\partial ^2}{\partial z_j\partial z_l}\log \det \left( g_{\alpha \overline{\beta }}\right) +q\frac{\partial f}{\partial z_l}\frac \partial {\partial
z_j}\log \det \left( g_{\alpha \overline{\beta }}\right) \right) }\det
\left( g_{\alpha \overline{\beta }}\right) ^q\nonumber\\
&&+g^{i\overline{j}}\left( \frac{\partial f}{\partial z_i}+qf\frac \partial
{\partial z_i}\log \det \left( g_{\alpha \overline{\beta }}\right) \right)
\times\nonumber\\ 
&&\quad\overline{\left( q\frac{\partial f}{\partial z_l}\frac{\partial ^2}{\partial 
\overline{z}_k\partial z_j}\log \det \left( g_{\alpha \overline{\beta }}\right) +qf\frac{\partial ^3}{\partial \overline{z}_k\partial z_l\partial
z_j}\log \det \left( g_{\alpha \overline{\beta }}\right) \right) }\det
\left( g_{\alpha \overline{\beta }}\right) ^q\ . \nonumber
\end{eqnarray}
Thus we obtain%
\begin{eqnarray}
\bigtriangleup \left\| \nabla s\right\| ^2&=&\left\| \nabla ^2s\right\| ^2+Ric\left(\nabla s,\nabla s\right)-trace\ C(L)\left\| \nabla s\right\| ^2+\left\| C(L)\right\| ^2\cdot \left\| s\right\| ^2\nonumber\\
&&+2Re\left\langle \nabla s,\nabla \left( trace\ C(L)\right) \otimes s\right\rangle -2C(L)\left(\nabla s,\nabla s\right)\ .\nonumber
\end{eqnarray}Suppose the sectional curvature is bounded by a positive
constant $K_0$. Thus the above Bochner formula simplifies to%
\begin{eqnarray}
\label{3.13}
\bigtriangleup \left\| \nabla s\right\| ^2&\geq&-CqK_0\left\| \nabla s\right\| ^2-Cq\left\| \nabla s\right\| \cdot \left\| s\right\| \cdot \left\| \nabla R\right\|\nonumber\\
&\geq&-C\left\| \nabla s\right\| ^2-\left\| s\right\| ^2\cdot \left\| \nabla R\right\| ^2\ .
\end{eqnarray}From Proposition 2.3 we know that $\left\| s\right\| $ is
bounded and $\left\| \nabla s\right\| $ is square--integrable. But the term
containing $\left\| \nabla R\right\| $ cause the difficulty to derive a
pointwise estimate for $\left\| \nabla s\right\| $ by using the mean value
inequality since at a priori we have no control on $\left\| \nabla R\right\| 
$. To overcome the difficulty we consider the Hamilton's Ricci flow equation
on $M,$%
\begin{equation}
\label{3.14}\left\{ 
\begin{array}{ll}
\bigbreak\displaystyle\frac{\partial g_{i\overline{j}}(x,t)}{\partial t}%
=-R_{i\overline{j}}(x,t)\ , & \qquad x\in M\ ,\quad t>0\ , \\ 
g_{i\overline{j}}(x,0)=g_{i\overline{j}}(x)\ , & \qquad x\in M\ . 
\end{array}
\right. 
\end{equation}
We may assume that the curvature tensor $R_m$ is also bounded by the
positive constant $K_0$. From \cite{Sh} we know that there exists a constant 
$\delta >0$, depending only on $K_0$ and the dimension $n$, such that the
Ricci flow equation (\ref{3.14}) has a solution $g_{i\overline{j}}(x,t)$ on $%
M\times [0,\delta ]$ with the following estimates 
\begin{equation}
\label{3.15}\frac 12g_{i\overline{j}}(\cdot )\leq g_{i\overline{j}}(\cdot
,t)\leq 2g_{i\overline{j}}(\cdot )\ , 
\end{equation}
\begin{equation}
\label{3.16}\left\| R_m(\cdot ,t)\right\| _t\leq 2K_0\ , 
\end{equation}
\begin{equation}
\label{3.17}\left\| \nabla ^tR_m(\cdot ,t)\right\| _t\leq \frac C{\sqrt{t}}\
, 
\end{equation}
where $\left\| \cdot \right\| _t$ and $\nabla ^t$ are the norm and the
covariant derivature with respect to the metric $g_{i\overline{j}}(\cdot ,t)$%
, the constant $C$ depending only on $K_0$ and the dimension.

Denote by $B_t(x_0,1)$ the geodesic ball of radius 1 and $Vol_t(B_t(x_0,1))$ its volume with respect to the metric $g_{i\overline{j}}(\cdot ,t)$.
Since the sectional curvature is positive and bounded at $t=0$, this implies
that the injectivity radius has a positive lower bound and then we have the
following estimate%
$$
Vol\left( B_0\left( x_0,\frac 1{\sqrt{2}}\right) \right) \geq \beta >0 
$$
for some positive constant depending only on $K_0$ and $n$. By using (\ref
{3.15}), it is easy to see%
\begin{eqnarray}
\label{3.18}
Vol_t(B_t(x_0,1))&\geq&\int_{B_0\left( x_0,\frac 1{\sqrt{2}}\right) }\det \left( g_{i\overline{j}}(\cdot ,t)\right)\nonumber\\
&\geq&\frac 1{2^n}Vol\left( B_0\left( x_0,\frac 1{\sqrt{2}}\right) \right)\nonumber\\
&\geq&\frac \beta {2^n}\ .
\end{eqnarray}Furthermore, we know from \cite{Sh} that the positivity of
holomorphic bisectional curvature is preserved under the Ricci flow (\ref
{3.14}). In particular, the Ricci curvature is positive for $t\in [0,\delta
].$

Note that $\left\| s\right\| $ is bounded and square--integrable. The estimate (\ref{3.15}) simply says that the metrics
$g_{i\overline j}(\cdot,t)$ are equivalent for $t\in[0,\delta]$. It follows that $s$ is still bounded and square--integrable
holomorphic section of $K^{-q}$ for $t\in [0,\delta ].$ Then Proposition 2.3 gives us
the $L^2$ gradient estimate 
\begin{equation}
\label{3.19}\left\| \nabla ^ts\right\| _{L^2\left( g_{i\overline{j}}(\cdot
,t)\right) }<+\infty \ . 
\end{equation}

On the other hand, by (\ref{3.13}), (\ref{3.16}), and (\ref{3.17}), we have%
$$
\bigtriangleup _\delta \left\| \nabla ^\delta s\right\| _\delta ^2\geq
-B \left\| \nabla ^\delta s\right\| _\delta ^2-B 
$$
for some positive constant $B$, where $\bigtriangleup _s$ is the Laplacian
operator of the metric $g_{i\overline{j}}(\cdot ,\delta )$. Consider $M$
equipped with the metric $g_{i\overline{j}}(\delta )$ and let $G=e^{\sqrt{B}%
\tau }\left( \left\| \nabla ^\delta s\right\| _\delta ^2+1\right) $ be a
function defined on the product manifold $\widetilde{M}=M\times {\R}$,
equipped with the product metric. Let $\widetilde{\bigtriangleup }%
=\bigtriangleup _\delta +\frac{\partial ^2}{\partial \tau ^2}$ be the
Laplacian operator. It is clear that $\widetilde{M}$ has nonnegative Ricci
curvature and we have%
$$
\widetilde{\bigtriangleup }G\geq 0\ . 
$$
From (\ref{3.18}) and (\ref{3.19}), by using the mean value inequality for
the subharmonic function $G$ on $\widetilde{M}$ exactly as in the proof of
Proposition 2.3, we see that 
\begin{equation}
\label{3.20}\sup \limits_{x\in M}\left\| \nabla ^\delta s\right\| _\delta
^2(x)<+\infty \ . 
\end{equation}

By (\ref{3.12}) we have%
\begin{eqnarray}
\nabla ^\delta s-\nabla s&=&qf\frac \partial {\partial z_i}\log \frac{\det \left( g_{\alpha \overline{\beta }}(\cdot ,\delta )\right) }{\det \left( g_{\alpha \overline{\beta }}(\cdot ,0)\right) }dz_i\otimes \left( \frac \partial {\partial z_1}\wedge \cdots \wedge \frac \partial {\partial z_n}\right) ^q\nonumber\\
&=&q\partial F\otimes s\ ,\nonumber
\end{eqnarray}where $F=\log \det \left( g_{\alpha \overline{\beta }}(\cdot
,\delta )\right) -\log \det \left( g_{\alpha \overline{\beta }}(\cdot
,0)\right) $. By using the Ricci flow equation (\ref{3.14}) one readily sees
that%
$$
F(x,t)=-\int_0^\delta R(x,t)dt\ . 
$$
Thus by combining (\ref{3.15}) and (\ref{3.17}), we have%
\begin{eqnarray}
\left\| \partial F\right\| _0&\leq&\sqrt{2}\int_0^\delta \left\| \nabla ^tR(\cdot ,t)\right\| _tdt\nonumber\\
&\leq&\sqrt{2}\int_0^\delta \frac C{\sqrt{t}}dt\nonumber\\
&\leq&const.\ ,\nonumber
\end{eqnarray}which implies that%
$$
\left\| \nabla ^\delta s-\nabla s\right\| _0\leq const.\left\| 
s\right\| _0\ , 
$$
and then%
\begin{eqnarray}
\left\| \nabla s\right\| _0&\leq&\left\| \nabla ^\delta s-\nabla s\right\| _0+\left\| \nabla ^\delta s\right\| _0\nonumber\\
&\leq&const.\left\| s\right\| _0+2^{\frac{nq+1}2}\left\| \nabla ^\delta s\right\| _\delta\nonumber\\
&<&+\infty\nonumber
\end{eqnarray}by (\ref{3.15}) and (\ref{3.20}). Therefore we have completed
the proof of the proposition.

\hfill$\Box $

We now can apply Proposition 3.3 to the smooth holomorphic curve $S$ which lies in
the intersection of the zero divisors of the $n-1$ holomorphic sections $%
t_1,\cdots ,t_{n-1}$, where $t_i\in \Gamma ^2(M,K^{-q}).$

By using Proposition 3.3 and (\ref{3.11}) we see that%
\begin{eqnarray}
\int\limits_{{\P}TM}\zeta _{\varepsilon _1}(\pi ^{*}t_1)\wedge \cdots \wedge \zeta
_{\varepsilon _{n-1}}(\pi ^{*}t_{n-1})\wedge\eta _{\varepsilon _n}(t_0\nabla t_1-t_1\nabla t_0)&&\nonumber\\
\wedge \cdots\wedge\eta _{\varepsilon _{2n-2}}(t_0\nabla t_{n-1}-t_{n-1}\nabla t_0)\wedge\nu&\leq& C\ q^{2n-2}\int_MRic^n, \nonumber
\end{eqnarray}
and by Fatou's Lemma and the Poincar\'e--Lelong equation, 
$$
\int_{\widehat{S}}\nu \leq Cq^{2n-2}\int_MRic^n<+\infty \ . 
$$

On the other hand, a calculation as in Mok \cite{Mo2} easily leads to the
following formula%
$$
\int_{\widehat{S}}\nu =-\int_SK(x)+2\int_SRic\ , 
$$
where $K(x)$ is the Gaussian curvature of $S$. Hence we have 
$$
\int_SK(x)\geq -Cq^{2n-2}\int_MRic^n\ , 
$$
where $C$ is a positive constant depending only on $n.$ Thus we prove the
following result.\vskip 3mm{\bf \underline{Proposition 3.5}}\quad Let $M$ be
assumed as in the Main Theorem. Suppose $S$ is a smooth holomorphic curve which
lies in the intersection of $n-1$ zero divisors of $t_1,\cdots ,t_{n-1}\in
\Gamma ^2(M,K^{-q})$, $(q\geq 2)$. Then the Gauss--Bonnet integral satisfies,%
\begin{equation}
\label{3.21}
\int_SK(x)\geq -Cq^{2n-2}\int_MRic^n\ , 
\end{equation}
where $C$ is a positive constant depending only on $n.$\hfill$\Box $

\section*{\S 4. Proof of the Main Theorem}

\setcounter{section}{4} \setcounter{equation}{0} \qquad This section is
devoted to the proof of the Main Theorem in the introduction. Let $M$ be
assumed as in the Main Theorem. From Section 2 we know that the space $%
\bigcup\limits_{q>0}\Gamma ^2\left( M,K^{-q}\right) $ forms a {{\Z}$^{+}$}%
--graded algebra, gives local holomorphic coordinates at any point in $M$
and separates points of $M$. Applying Fatou's lemma to Proposition 3.2 with $%
k=1$, we have 
\begin{equation}
\label{4.1}\int_{[s=0]}Ric^{n-1}\leq \lim \limits_{\ \varepsilon \rightarrow } 
\negthinspace \inf _{0\ \ }\int_M\zeta _\varepsilon (s)\wedge Ric^{n-1}\leq
Cq\int_MRic^n 
\end{equation}
for any $s\in \Gamma ^2\left( M,K^{-q}\right) $, where $C$ is a positive
constant depending only on $n$, and $[s=0]$ denotes the zero divisor of $s$
counting multiplicity. Fix a point $x_0\in M$. Since the Ricci curvature is
positive, it follows from (\ref{4.1}) and the inequality of Bishop--Lelong
that the multiplicity estimate%
$$
mult\left( \left[ s=0\right] ,x_0\right) \leq C_1q 
$$
holds for all $s\in \Gamma ^2\left( M,K^{-q}\right) $, where the positive
constant $C_1$ may depend on $x_0$ but independent of $s$ and $q$. Then it
follows from a standard argument by considering the local Taylor expansion
of $s\in \Gamma ^2\left( M,K^{-q}\right) $ at the fixed point $x_0$ ( c.f.
the proof of Proposition 5.1 in \cite{Mo1} ) that 
\begin{equation}
\label{4.2}\dim_{{\C}}\Gamma ^2\left( M,K^{-q}\right) \leq C_2q^n 
\end{equation}
where $C_2$ is independent of $q$.

Denote by $R\left( M,K^{-1}\right) $ the subset of meromorphic functions on $%
M$ obtained by taking quotients in $\bigcup\limits_{q>0}\Gamma ^2\left(
M,K^{-q}\right) $ of the same degree, i.e., 
$$
R\left( M,K^{-1}\right) =\left\{ \left. \ \frac st\ \right| \ s,t\in \Gamma
^2\left( M,K^{-q}\right) ,\ where\ q\in {\Z}^{+}\ and\ t\not \equiv 0\
\right\} \ . 
$$
Since $\bigcup\limits_{q>0}\Gamma ^2\left( M,K^{-q}\right) $ forms a graded
algebra under addition and multiplication, it thus is easy to see that $%
R\left( M,K^{-1}\right) $ forms a subfield of the field of meromorphic
functions on $M$. Let $s_0,s_1,\cdots ,s_n$ be $n+1$ holomorphic sections in 
$\Gamma ^2\left( M,K^{-q}\right) $ ( for some $q>0$ large enough )
constructed in Section 2 such that $s_1/s_0,\cdots ,s_n/s_0$ give a local
holomorphic coordinate system at $x_0$. By the estimate (\ref{4.2}), the
classical argument of Poincar\'e--Siegel ( c.f. e.g. Mok \cite{Mo1} ) shows
that $R\left( M,K^{-1}\right) $ is a finite extension field over ${\C}\left( 
\frac{s_1}{s_0},\cdots ,\frac{s_n}{s_0}\right)$. And by the primitive
element theorem, the subfield $R\left( M,K^{-1}\right) $ is given by%
$$
R\left( M,K^{-1}\right) ={\C}\left( \frac{s_1}{s_0},\cdots ,\frac{s_n}{s_0}%
,g\right)\ , 
$$
where $g$ is some meromorphic function in $R\left( M,K^{-1}\right) $ and is
algebraic over ${\C}\left( \frac{s_1}{s_0},\cdots ,\frac{s_n}{s_0}\right)$.
Also by taking common denominators we may assume $g=s_{n+1}/s_0$, where $%
s_0,s_1,\cdots ,s_n,s_{n+1}\in \Gamma ^2\left( M,K^{-q}\right) $ for some $%
q>0$ large enough. Now consider the mapping $F:M\rightarrow {\P}^{n+1}$
defined by%
$$
F(x)=\left[ s_0(x),s_1(x),\cdots ,s_n(x)\right] ,\qquad for\quad x\in M\ . 
$$
Since $g$ is algebraic, the minimal equation satisfied by $g$ over ${\C}
\left( \frac{s_1}{s_0},\cdots ,\frac{s_n}{s_0}\right)$ can be given by%
$$
g^p+\sum\limits_{0\leq j\leq p-1}R{_j\left( \frac{s_1}{s_0},\cdots ,\frac{%
s_n }{s_0}\right) g^j=0\ ,} 
$$
where $R_j,\ 1\leq j\leq p-1$, are rational functions of $n$ variables.
After clearing denominators, we see that $s_0,s_1,\cdots ,s_{n+1}$ satisfy a
homogeneous equation%
$$
P\left( s_0,s_1,\cdots ,s_{n+1}\right) =0\ . 
$$
Let $Z_0$ be the hypersurface of ${\P}^{n+1}$ defined by%
$$
Z_0=\left\{ \left. \ \left[ s_0,s_1,\cdots ,s_{n+1}\right] \in {\P}^{n+1}\
\right| \ P\left( s_0,s_1,\cdots ,s_{n+1}\right) =0\ \right\} \ , 
$$
and let $Z$ be the connected component of $Z_0$ containing $F\left(
M\backslash Q\right) $ where $Q$ is the set of indeterminancy of $F.$

In the following we will show that $F$ is an ``almost injective'' and
``almost surjective'' map to $Z$ and we can desingularize $F$ to obtain a
biholomorphic map from $M$ onto a quasi--projective variety by adjoining a
finite number of holomorphic plurianticanonical sections.

First of all, we claim that $Z$ is irreducible and $F$ is ``almost
injective'', i.e., there exists a subvariety $V$ of $M$ such that $%
F|_{M\backslash V}:M\backslash V\rightarrow Z$ is an injective locally
biholomorphic mapping. Indeed, take $V$ to be the union of the branching
locus and the base locus of $F$ and $F^{-1}(Sing(Z))$, here $Sing(Z)$
denotes the singular set of $Z$. It is clear that $F$ is locally
biholomorphism on $M\backslash V$. That $F$ is also injective there follows
from the fact that $\bigcup\limits_{q>0}\Gamma ^2\left( M,K^{-q}\right) $
separates points and $s_1/s_0,\cdots ,s_n/s_0,s_{n+1}/s_0$ generate $R\left(
M,K^{-1}\right) $. To see the irreducibility of $Z $, note that $M\backslash
(Q\cup (SingZ))$ is connected since $Q\cup (SingZ)$ has real codimension 2.
Here $\overline{F\left( M\backslash (Q\cup (SingZ))\right) }$ is irreducible
( as its set of smooth points is connected ). Since $F\left( M\backslash
Q\right) \subset \overline{F\left( M\backslash (Q\cup (SingZ))\right) }$, by
the definition of $Z$, it must be irreducible.

We now make a remark for the subvariety $Z$. We have seen that $F$ is a
biholomorphic embedding from $M\backslash V$ into its image $F\left(
M\backslash V\right) $ in $Z$. As in Mok \cite{Mo2}, $F$ is a birational
mapping from $M$ to $Z$ with respect to the subfield $R\left(
M,K^{-1}\right) $, i.e., the rational function field over $Z$ is isomorphic
to $R\left( M,K^{-1}\right) $ via the pull back $F^{*}$. Therefore for any
projective desingularization $\pi :Z^{\prime }\rightarrow Z$, where $%
Z^{\prime }$ is a smooth projective variety, $\pi ^{-1}\circ F$ is still a
birational mapping from $M$ to $Z^{\prime }$ with respect to $R\left(
M,K^{-1}\right) $. Since $\pi ^{-1}$ is a birational map ( in the usual
sense ), the composition $\pi ^{-1}\circ F$ is easily seen to be defined by
sections in $\Gamma ^2\left( M,K^{-q}\right) $ for some large positive
integer. Henceforth at each step when we desingularize $F$ by adding
plurianticanonical sections, we may assume $Z$ is smooth.

Next, we come to the ``almost surjectivity'' of $F$, i.e., the image $%
F\left( M\backslash V\right) $ contains a Zariski--open subset of $Z$. Since
the sectional curvature of $M$ is positive, by the result of Greene and Wu 
\cite{GW1}, $M$ is Stein. Thus as in Mok--Zhong \cite{MZ} ( see also To \cite
{T} ) by using Simha's criterion \cite{Sim} for a domain of holomorphy
contained in the unit polydisc to be Zariski open, the proof of the ``almost
surjectivity'' of $F$ is reduced to show that the cardinality of the
intersections of the complement of $F\left( M\backslash V\right) $ with the
generic curves obtained by intersecting $Z$ with $(n-1)$ hyperplane sections
of ${\P}^{n+1}$ are bounded. The following argument is basically due to
Mok--Zhong ( see Proposition (2.2.2) and Proposition (2.2.3) in \cite{MZ} ).

Let $S_0\subset Z$ be a curve obtained by intersecting $Z$ with $(n-1)$
hyperplane sections, which correspond to $(n-1)$ holomorphic sections, say $%
t_1,\cdots ,t_{n-1}$, in $\Gamma ^2\left( M,K^{-q}\right) $. Assume that the
hyperplane sections and $Z$ intersect at normal crossings at $Q\in Z$ and
denote by $S$ the irreducible component of $S_0$ containing $Q$. Let $D$
denote $F\left( M\backslash V\right) \cap S$ and assume that $D$ is nonempty
for the generic $S$. Define $\Sigma $ to be the closure of $F^{-1}(D)$ on $M$
and $f$ to be the meromorphic extension of $F|_{F^{-1}(D)}$ to $\Sigma $. In
order to prove that the cardinality of $S\backslash D$ is bounded, we
proceed to prove that $f(\Sigma )\backslash D$ and $S\backslash f(\Sigma )$
are finite set of points, giving estimates of cardinality of these two sets
at the same time. Note that $\Sigma \backslash f^{-1}(D)$ lies in the
intersection of the zero divisors of $t_1,\cdots ,t_{n-1}$ and $V$, and $V$
is clearly contained in a divisor determined by the sections $s_0,s_1,\cdots
,s_n$ which define the mapping $F$. Applying Fatou's lemma to the B\'ezout
estimate (\ref{3.1}) in Proposition 3.2 we see that $\Sigma \backslash
f^{-1}(D)$ consists of a finite number of points and the bound on the
cardinality is independent of $\Sigma $. This gives an upper bound for $%
f(\Sigma )\backslash D$. To get an upper bound for $S\backslash f(\Sigma )$
one can pass to a normalization $\sigma :S^{\prime }\rightarrow S$ and prove
the corresponding statement for $S^{\prime }$ and the Riemann surface $%
D^{\prime }=\sigma ^{-1}(f(\Sigma ))$, which is equipped with the induced (
possible degenerate ) Hermitian metric $\left( f^{-1}\circ \sigma \right)
^{*}\left( \omega |_\Sigma \right) $. The degeneracies occur at some points
corresponding to singularities of $\Sigma $. Such singularities can either
come from singularities of $D$ itself or from the intersection of $\Sigma $
with $V$. The first set is finite since $Z $ is algebraic, and the second
set is also finite by the B\'ezout estimate (\ref{3.1}) in Proposition 3.2.
Thus $\left( f^{-1}\circ \sigma \right) ^{*}\left( \omega |_\Sigma \right) $
is degenerate at at most a finite set of points $\left\{ p_i\right\}
_{i=1}^r $. Modify the metric in small compact neighborhoods of these points
to get a smooth metric $\mu $. Then the Gauss--Bonnet integrals of these two
metrics are related by%
$$
\int_{D^{\prime }}c_1(D^{\prime },\mu )=\int_{D^{\prime }\backslash
\{p_i\}}c_1\left( D^{\prime },\left( f^{-1}\circ \sigma \right) ^{*}\left(
\omega |_\Sigma \right) \right) -r\ ,
$$
where $c_{1}$ is denoted by the associated Chern form.
Recall that the estimate (\ref{3.21}) in Proposition 3.5 says%
$$
\int_{D^{\prime }\backslash \{p_i\}}c_1\left( D^{\prime },\left( f^{-1}\circ
\sigma \right) ^{*}\left( \omega |_\Sigma \right) \right) \geq
-Cq^{2n-2}\int_MRic^n\ , 
$$
which implies that%
$$
\int_{D^{\prime }}c_1(D^{\prime },\mu )>-C 
$$
for a constant $C$ independent of the choice of hyperplane sections. By the
well--known result of Huber \cite{Hu} we have%
$$
2genus(S^{\prime })+card\left( S^{\prime }\backslash D^{\prime }\right)
-2\leq -\frac 1{2\pi }\int_{D^{\prime }}c_1(D^{\prime },\mu )<C\ . 
$$
This gives the upper bound for the number of points in $S^{\prime
}\backslash D^{\prime }$ and hence of $S\backslash f(\Sigma )$. Thus we
obtain the uniform upper bound for $S\backslash D$. Therefore $F(M\backslash
V)$ is Zariski open in $Z.$

Now we are in position to desingularize $F$ to a biholomorphism from $M$ to
a quasi--projective variety. Recall that for any fixed point $x_0\in V$, one
can choose $n+1$ holomorphic plurianticanonical sections $s_0^{\prime
},s_1^{\prime },\cdots ,s_n^{\prime }\in \Gamma ^2\left( M,K^{-q^{\prime
}}\right) $ which give the holomorphic coordinate at $x_0$. Although the
integer $q^{\prime }$ may be different from the integer $q$ in the mapping $%
F $, we can add the sections $s_0^{\prime },s_1^{\prime },\cdots
,s_n^{\prime } $ to $F$ to get a birational map $F^{\prime }$ by composing
the sections and $F$ with a Veronese map as done in \cite{MZ}, \cite{T}.
Denote by $V^{\prime }$ the union of the branching locus and the base locus
of $F^{\prime }$. Then $V^{\prime }\stackrel{\subset }{\neq }V$. If for each
such a birational map, the union of the branching locus and the base locus
has only finite number of irreducible components of maximal dimension, then
it is clear that the birational map $F$ will become a biholomorphism after
adding a finite number of holomorphic plurianticanonical sections in $%
\bigcup\limits_{q>0}\Gamma ^2\left( M,K^{-q}\right) .$

Recall a result of Demailly \cite{D} which states that if $\dim V \leq p$, then
the relative cohomology group $H^{2(n-p)}(M,M\backslash V,{\R)}$ is
isomorphic to ${\R}^J$, where $J$ is the number of irreducible components of
dimension $p$ in $V$. As a consequence, if $H^{2(n-p)}(M,{\R)}$ and $%
H^{2(n-p)-1}(M\backslash V,{\R)}$ are finite dimensional, then the number of
irreducible branches of $V$ of complex dimension $p$ is finite. In our case, 
$M$ is diffeomorphic to ${\R}^{2n}$ and $M\backslash V$ is Zariski open in
an algebraic variety. So $M\backslash V$ is of finite topological type.
Hence we can desingularize the map $F$ into a biholomorphism from $M$ to a
quasi--projective variety.This proves the first assertion of the Main Theorem. 

To prove the second
assertion of the Main Theorem we recall a theorem of Ramanujam \cite{R}
which states that a quasi--projective variety homeomorphic to ${\R}%
^4$ is biregular to ${\C}^2$. Thus $M$ is biholomorphic to ${\C}^2$ when the
complex dimension $n=2.$

Therefore the proof of the Main Theorem is completed.\hfill$\Box $

\end{document}